\documentclass{ifacconf}
\usepackage{graphicx} 
\usepackage
{
    amssymb,
    amsmath,
    dsfont, 
    xcolor,
    mathtools,
    nicefrac,
    array,
    enumerate,
    siunitx,
    tikz,
}
\usepackage{subcaption}
\usepackage{natbib}

\usetikzlibrary{shapes,arrows}
\usetikzlibrary{calc} 

\usepackage[utf8]{inputenc}

\graphicspath{{resources/}}

\makeatletter
\newcommand\RedeclareMathOperator{%
  \@ifstar{\def\rmo@s{m}\rmo@redeclare}{\def\rmo@s{o}\rmo@redeclare}%
}
\newcommand\rmo@redeclare[2]{%
  \begingroup \escapechar\m@ne\xdef\@gtempa{{\string#1}}\endgroup
  \expandafter\@ifundefined\@gtempa
     {\@latex@error{\noexpand#1undefined}\@ehc}%
     \relax
  \expandafter\rmo@declmathop\rmo@s{#1}{#2}}
\newcommand\rmo@declmathop[3]{%
  \DeclareRobustCommand{#2}{\qopname\newmcodes@#1{#3}}%
}
\@onlypreamble\RedeclareMathOperator
\makeatother

\newcommand{\N}{\mathds{N}}
\newcommand{\R}{\mathds{R}}
\newcommand{\Rp}{\R_{\geq0}}

\newcommand{\Impl}{\Longrightarrow}

\newcommand{\fa}{\ \forall \, }
\newcommand{\ex}{\ \exists \, }

\newcommand{\rbl}{\left (}
\newcommand{\rbr}{\right )}

\newcommand{\al}{\left \langle}
\newcommand{\ar}{\right \rangle}

\newcommand{\nl}{\left\|}
\newcommand{\nr}{\right\|}
\newcommand{\cbl}{\left\lbrace }
\newcommand{\cbr}{\right\rbrace }

\newcommand{\Norm}[2][ ]{\nl #2 \nr_{#1}}
\newcommand{\SNorm}[1]{\Norm[\infty]{#1}}

\newcommand{\setdef}[2]{\cbl\ #1\ \left|\ \vphantom{#1} #2\ \right.\cbr}

\newcommand{\cD}{\mathcal{D}}

\newcommand{\cF}{\mathcal{F}}
\newcommand{\cG}{\mathcal{G}}
\newcommand{\cN}{\mathcal{N}}

\newcommand{\cT}{\mathcal{T}}

\DeclareMathOperator*{\rf}{ref}
\DeclareMathOperator*{\esssup}{ess\,sup}

\DeclareMathOperator*{\loc}{loc}

\newcommand{\eps}{\varepsilon}

\newcommand{\me}{\mathrm{e}}

\newcommand{\cC}{\mathcal{C}}
\newcommand{\cE}{\mathcal{E}}
\newcommand{\oT}{\mathbf{T}}
\DeclareMathOperator*{\graph}{graph}

\RedeclareMathOperator*{\Im}{Im}
\RedeclareMathOperator*{\Re}{Re}
\renewcommand{\phi}{\varphi}

\newcommand{\dd}[2][ ]{\tfrac{\text{\normalfont d}#1}{\text{\normalfont d}#2}}

\makeatletter
\@ifundefined{ifCustomTheorems}{}
{
    \setcounter{section}{0}
    \counterwithout{equation}{section}
    \counterwithout{figure}{section}
    \newtheorem{definition}{Definition}[section]
    \theoremstyle{definition}
    \newtheorem{remark}[definition]{Remark}\Crefname{remark}{Remark}{Remarks}
    \Crefname{algo}{Algorithm}{Algorithms}
    \Crefname{example}{Example}{Examples}
    \theoremstyle{plain}
    \Crefname{prop}{Proposition}{Propositions}
    \Crefname{corollary}{Corollary}{Corollaries}
    \Crefname{assertion}{Assertion}{Assertions}
    \newtheorem{theorem}[definition]{Theorem}\Crefname{theorem}{Theorem}{Theorems}
    \newtheorem{lemma}[definition]{Lemma}\Crefname{lemma}{Lemma}{Lemmata}
}
\makeatother

\newtheorem{theorem}{Theorem}[section]
\newtheorem{definition}{Definition}[section]
\newtheorem{remark}{Remark}[section]
\newtheorem{lemma}{Lemma}[section]

\DeclareOldFontCommand{\rm}{\normalfont\rmfamily}{\mathrm}

\makeatletter
\renewcommand*\env@matrix[1][*\c@MaxMatrixCols c]{%
  \hskip -\arraycolsep
  \let\@ifnextchar\new@ifnextchar
  \array{#1}}
\makeatother

\newcommand{\OpChi}[1][r]{\chi_{#1}}
\newcommand\xqed{\leavevmode\unskip\penalty9999 \hbox{}\nobreak\hfill \quad\hbox{\qed}}

\begin{document}
\begin{frontmatter}

\title{A low-complexity funnel control approach for non-linear systems of higher-order}

\thanks[footnoteinfo]{This work was supported by the German Research Foundation - Deutsche Forschungsgesellschaft 
(DFG; project number 471539468)}

\author[First]{Dario Dennstädt} 

\address[First]{Institut f\"ur Mathematik, Universit\"at Paderborn,\\ Warburger Stra{\ss}e~100, 33098~Paderborn, Germany\\
   (e-mail: {dario.dennstaedt@uni-paderborn.de})}

\begin{abstract}       
We address the problem of output reference tracking for unknown non-linear
multi-input, multi-output systems described by functional differential
equations. This class of systems includes those with a strict relative degree,
and bounded-input bounded-output (BIBO) stable internal dynamics. The objective
is to ensure that the tracking error evolves within a prescribed performance
funnel. To achieve this, we propose a novel model-free adaptive controller with
lower complexity than existing funnel control methods, avoiding the use of
non-linearities in intermediate error signals. We establish the feasibility and
effectiveness of the proposed design through theoretical analysis and
demonstrate its performance with simulations, comparing it to previous
approaches.
\end{abstract}


\end{frontmatter}

\section{Introduction}
\emph{Funnel control} is a model-free, high-gain output feedback technique first
introduced in~\cite{IlchRyan02b}. It ensures output tracking of reference
signals within prescribed performance bounds for a broad class of non-linear
multi-input, multi-output systems. The approach solely invokes certain
structural assumptions about the system, namely stable internal dynamics and a
known relative degree with a sign-definite high-frequency matrix. Under these
conditions, the adaptive controller offers robustness against disturbances and
guarantees specified transient behaviour without relying on explicit knowledge
about system to be controlled. Funnel control therefore proved useful for
tracking problems in various applications such as DC-link power flow control,
see~\cite{SenfPaug14}, control of industrial servo-systems, see~\cite{Hackl17},
and temperature control of chemical reactor models, see~\cite{IlchTren04}.

Recent research on funnel control has focused on addressing practical challenges
that might arise in real-world applications. These include issues like
measurement losses, see \cite{BergerMeasurement}, input constraints, see
\cite{BergerInput}, reliance on sampled measurement data, see \cite{Lanza24},
and leveraging existing model knowledge, see \cite{BergerFMPC}. For a
comprehensive overview of the literature, we also recommend the survey
paper~\cite{berger2023funnel}. Despite these advancements, a persistent
challenge for funnel control and, more general, for high-gain adaptive control
has been the handling of systems with relative degree larger than one, for the
latter see e.g.~\cite{MORSE96}.

To address the problem of output tracking within prescribed boundaries for
higher-order systems, \cite{BergLe18a} proposed a control strategy utilising
auxiliary error variables, with the number of variables corresponding to the
system's relative degree. Each auxiliary error signal is recursively defined and
includes a reciprocal penalty term that depends non-linearly on the preceding
variable. These signals evolve within their respective funnel boundaries,
maintaining the desired tracking performance. However, due to their reciprocal
nature, error gains can become very large when an auxiliary variable nears its
boundary, posing potential numerical stability issues. Although the controller
design was further simplified in~\cite{BergIlch21}, the general approach remains
similar, and the controllers lack the low complexity of the original design for
systems of order one.

\emph{Prescribed performance control} is a relative to funnel control and was
first  proposed in \cite{BECHLIOULIS2008}. The tracking within prescribed
boundaries is achieved by transforming the system via so-called performance
functions. While the original controller design incorporates neural networks to
approximate unknown non-linearities, \cite{BECHLIOULIS2014} proposed an
approximation-free control scheme for systems in pure feedback form. Similar to
funnel control, recursively defined intermediate control signals are used,
matching the number of the system's relative degree and depending non-linearly
on their predecessors.

In this paper, we propose a funnel controller of lower complexity for systems of
higher-order. While still relying on recursively defined auxiliary error
variables to handle the system's higher-order nature, the proposed approach
eliminates the use of time-varying reciprocal penalty terms, replacing them with
constant gains. This modification simplifies the controller design and has the
potential to mitigate numerical issues and enhance its practicality for
real-world applications.

\subsubsection{Nomenclature:}
$\N$ and $\R$ denote natural and real numbers, respectively.
$\N_0:=\N\cup\{0\}$ and $\Rp:=[0,\infty)$.
$\Norm{x}:=\sqrt{\al x,x\ar}$~denotes the Euclidean norm of $x\in\R^n$.
$\cC^p(V,\R^n)$ is the linear space of $p$-times continuously  differentiable
functions $f:V\to\R^n$, where $V\subset\R^m$, and $p\in\N_0\cup \{\infty\}$.
$\cC(V,\R^n):=\cC^0(V,\R^n)$.
On an interval $I\subset\R$,  $L^\infty(I,\R^n)$ denotes the space of measurable and essentially bounded
functions $f: I\to\R^n$ with norm $\SNorm{f}:=\esssup_{t\in I}\Norm{f(t)}$,
$L^\infty_{\text{loc}}(I,\R^n)$ the set of measurable and locally essentially bounded functions.
Furthermore, $W^{k,\infty}(I,\R^n)$ is the Sobolev space of all $k$-times weakly differentiable functions
$f:I\to\R^n$ such that $f,\dots, f^{(k)}\in L^{\infty}(I,\R^n)$.

\section{Problem formulation}
We consider non-linear multi-input multi-output systems with constant degree $r\in\N$ of the form 
\begin{equation} \label{eq:Sys}
    \begin{aligned}
    y^{(r)}(t) &= f \big(d(t),\oT(y,\ldots,y^{(r-1)} )(t), u(t)\big) \\
   y|_{[0,t_0]} &= y^0  \in \cC^{r-1}([0,t_0],\R^m)
    \end{aligned}
\end{equation}
 with $t_0\ge 0$, initial trajectory $y^0$, 
control input~$u\in L^\infty_{\loc}([t_0,\infty), \R^m)$, 
and output $y(t)\in\R^m$ at time $t\geq 0$.
In the case of $t_0=0$, we identify $\cC^{r-1}([0,t_0],\R^m)$
with the vector space $\R^{rm}$. 
Then, the initial value condition in~\eqref{eq:Sys} is replaced by
$(y(t_0),\dot{y}(t_0),\ldots, y^{(r-1)}(t_0))=y^0\in\R^{rm}$.
Note that $u$ and $y$ have the same dimension~$m\in\N$.
The system consists of the \emph{unknown} non-linear function 
$f\in\cC(\R^p\times \R^q \times \R^m,\R^m)$,
\emph{unknown} non-linear operator~$\oT:\cC(\Rp,\R^m)\to L^\infty_{\loc}([t_0,\infty),\R^q)$, 
and \emph{unknown} bounded disturbances $d\in L^\infty([t_0,\infty),\R^p)$.
The operator~$\oT$ is causal, locally Lipschitz and satisfies a bounded-input bounded-output property. 
It is characterised in detail in the following definition.
\begin{definition} \label{Def:OperatorClass} 
For $n,q\in\N$, and $t_0\geq 0$, the set $\cT^{n,q}_{t_0}$ denotes the class of operators $\oT:
\cC(\Rp,\R^n) \to L^\infty_{\loc} ([t_0,\infty), \R^{q})$
for which the following properties hold:
\begin{itemize}
    \item\emph{Causality}:  $\fa y_1,y_2\in\cC(\Rp,\R^n)$  $\fa t\geq t_0$:
    \[
        y_1\vert_{[0,t]} = y_2\vert_{[0,t]}
        \quad \Impl\quad
        \oT(y_1)\vert_{[t_0,t]}=\oT(y_2)\vert_{[t_0,t]}.
    \]
    \item\emph{Local Lipschitz}: 
    $\fa t \ge t_0 $ $\fa y \in \cC([0,t], \R^n)$ 
    $\ex \Delta, \delta, c > 0$ 
    $\fa y_1, y_2 \in \cC(\Rp, \R^n)$ with
    $y_1|_{[0,t]} = y_2|_{[0,t]} = y $ 
    and $\Norm{y_1(s) - y(t)} < \delta$,  $\Norm{y_2(s) - y(t)} < \delta $ for all $s \in [t,t+\Delta]$:
    \[
     \hspace*{-2mm}   \esssup_{\mathclap{s \in [t,t+\Delta]}}  \Norm{\oT(y_1)(s) \!-\! \oT(y_2)(s) }  
        \!\le\! c \ \sup_{\mathclap{s \in [t,t+\Delta]}}\ \Norm{y_1(s)\!-\! y_2(s)}\!.
    \] 
    \item\emph{Bounded-input bounded-output (BIBO)}:
    $\fa c_0 > 0$ $\ex c_1>0$  $\fa y \in \cC(\Rp, \R^n)$:
    \[
    \sup_{t \in \Rp} \Norm{y(t)} \le c_0 \ 
    \Impl \ \sup_{t \in [t_0,\infty)} \Norm{\textbf{T}(y)(t)}  \le c_1.
    \]
\end{itemize}
\end{definition}
Note that using this operator many physical phenomena such as \emph{backlash},
and \emph{relay hysteresis}, and \emph{non-linear time delays}
can be modelled by $\oT$ where $t_0$ corresponds to the initial delay, cf.~\cite[Sec.~1.2]{BergIlch21}.
In addition, systems with infinite-dimensional internal dynamics can be represented by~\eqref{eq:Sys}, see~\cite{BergPuch20}. 

While unknown, we assume the function $f$ in the system~\eqref{eq:Sys} to have the so-called\
\emph{high-gain} property~\cite[Definition 1.2]{BergIlch21}.
It is an essential assumption in high-gain adaptive control which, loosely speaking, 
ensures that, if a large enough control input is applied, then the system reacts sufficiently fast.
We summarise our assumptions and define the general system class under consideration
which is only parametrised by two natural numbers $m$ and $r$.
\begin{definition} \label{Def:system-class}
    We say that the system~\eqref{eq:Sys} belongs to the system class
    $\cN^{m,r}$, written $(d, f, \oT) \in\cN^{m,r}$, if, for $t_0\geq0$ and
    some $q\in\N$, the following holds: $d\in L^\infty([t_0,\infty),\R^p)$,
    ${\oT\in\cT^{rm,q}_{t_0}}$, and 
    $f\in\cC(\R^p\times \R^q \times \R^m,\R^m)$ has the \emph{high-gain property}, i.e.
    there exists~$\nu\in(0,1)$ such that,
    for every compact sets $K_p\subset \R^p$, $K_q\subset\R^q$,
    the function $\mathfrak{h}:\R\to\R$ defined by
    \[
    \mathfrak{h}(s):=\min\setdef{\langle v, f(\delta,z, -s v)\rangle \!\!}
    {\!\!\! 
    \begin{array}{ll}
         \delta\in K_p,  z \in  K_q,   \\
    v\in\R^m,\nu \leq \|v\| \leq 1
    \end{array}
    \!\!\!}
    \]
    satisfies $\sup_{s\in\R} \mathfrak{h}(s)=\infty$.
\end{definition}

\begin{remark}
    Under assumptions provided in~\cite[Cor.~5.6]{ByrnIsid91a},
    a non-linear system of the form
    \begin{equation}\label{eq:SysWithState}
    \begin{aligned}
       \dot{x}(t)  & = \tilde f(x(t)) + \tilde g(x(t)) u(t),\quad x(t_0)=x^0\in\R^n,\\
       y(t)        & = \tilde{h}(x(t)),
    \end{aligned} 
    \end{equation}
    with  non-linear functions~$\tilde f:\R^n\to \R^n$,
    $\tilde g:\R^n\to \R^{n\times m}$ and $\tilde{h} : \R^n \to \R^m$, 
    can be put in the form~\eqref{eq:Sys} via a coordinate transformation induced by a diffeomorphism~$\Phi:\R^n\to\R^n$.
    Then, the operator $\oT$ is the solution operator of the internal dynamics of the transformed system. 
    As we do not assume knowledge the function $f$ and the operator $\oT$ of the
    system~\eqref{eq:Sys}, all results presented in this paper can be directly
    applied to the system~\eqref{eq:SysWithState} without computing~$\Phi$.
\end{remark}

\subsection{Control objective}\label{Ssec:ContrObj}
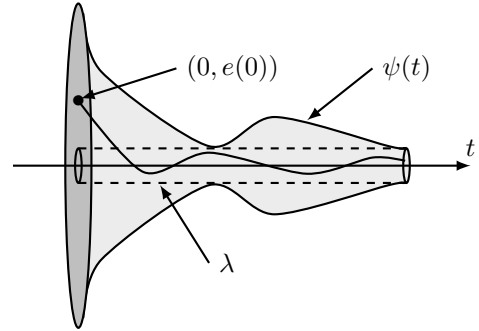
\begin{figure}[ht]
\hspace{1cm}
\begin{tikzpicture}[scale=0.43]
\tikzset{>=latex}
  \filldraw[color=gray!15] plot[smooth] coordinates {(0.15,4.7)(0.7,2.9)(4,0.6)(6,1.5)(9.5,0.6)(10,0.533)(10.01,0.531)(10.041,0.5) (10.041,-0.5)(10.01,-0.531)(10,-0.533)(9.5,-0.6)(6,-1.5)(4,-0.6)(0.7,-2.9)(0.15,-4.7)};
  \draw[thick] plot[smooth] coordinates {(0.15,4.7)(0.7,2.9)(4,0.6)(6,1.5)(9.5,0.6)(10,0.533)(10.01,0.531)(10.041,0.5)};
  \draw[thick] plot[smooth] coordinates {(10.041,-0.5)(10.01,-0.531)(10,-0.533)(9.5,-0.6)(6,-1.5)(4,-0.6)(0.7,-2.9)(0.15,-4.7)};
  \draw[thick,fill=lightgray] (0,0) ellipse (0.4 and 5);
  \draw[thick] (0,0) ellipse (0.1 and 0.533);
  \draw[thick,fill=gray!15] (10.041,0) ellipse (0.1 and 0.533);
  \draw[thick] plot[smooth] coordinates {(0,2)(2,-0.2)(4,0.4)(7,-0.25)(9,0.25)(10,0.15)};
  \draw[thick,->] (-2,0)--(12,0) node[right,above]{\normalsize$t$};
  \draw[thick,dashed](0,0.533)--(10,0.533);
  \draw[thick,dashed](0,-0.533)--(10,-0.533);
  \node [black] at (0,2) {\textbullet};
  \draw[->,thick](4,-3)node[right]{\normalsize$\lambda$}--(2.5,-0.6);
  \draw[->,thick](3,3)node[right]{\normalsize$(0,e(0))$}--(0.07,2.07);
  \draw[->,thick](9,3)node[right]{\normalsize$\psi(t)$}--(7,1.4);
\end{tikzpicture}
\caption{Error evolution in a funnel $\mathcal F_{\psi}$ with boundary $\psi(t)$. The figure is based on~\cite[Fig.~1]{BergLe18a}, edited for present purpose.}
\label{Fig:funnel}
\end{figure}

The objective is to design a output derivative feedback 
which achieves that, for any reference signal $y_{\rm ref}\in W^{r,\infty}(\Rp,\R^m)$, 
the output tracking error $e(t)=y(t)-y_{\rm ref}(t)$ evolves within a prescribed performance funnel
\begin{equation*}
    \mathcal{F}_{\psi} := \setdef{(t,e)\in\R_{\ge 0} \times\R^m} {\|e\| < \psi(t)}.
\end{equation*}
This funnel is determined by the choice of the function~$\psi$ belonging to
\begin{align*}
    \cG:=\setdef
        {\psi\in W^{1,\infty}(\Rp,\R)}
        {
          \begin{array}{l} \inf_{t\ge 0}  \psi(t) > 0,\\
          \exists\, \alpha, \beta>0\ \forall\, t\ge 0:\\
          \dot \psi(t) \ge -\alpha \psi(t) + \beta 
          \end{array}
        },
\end{align*}
see also Figure~\ref{Fig:funnel}.
Note that the tracking error evolving in $\mathcal{F}_{\psi}$ is not forced to asymptotically converge to zero. 
The specific application usually dictates the constraints on the tracking error and thus
indicates suitable choices for~$\psi$.

\section{Controller design}\label{Sec:ContrDes}
In order to achieve the control objective described in Section~\ref{Ssec:ContrObj}, 
we introduce auxiliary error variables. 
Define for $z=(z_1,\ldots,z_r)\in\R^{rm}$ and parameter $k\in\Rp$, the functions
$\xi_i:\R^{rm}\to\R^m$ recursively by 
\begin{equation}\label{eq:ErrorVars}
\boxed{
\begin{aligned}
    \xi_1(z) &:= z_1,\\
    \xi_{i+1}(z) &:= \xi_{i}(\sigma(z))+k\xi_{i}(z)\\
\end{aligned}
}
\end{equation}
for $i=1,\ldots, r-1$, where
\[
    \sigma:\R^{rm}\to\R^{rm},\ \sigma(z_1,\ldots,z_r)=(z_2,\ldots,z_r,0)
\]
is the left shift operator.
\begin{remark}
We utilise the notation 
\begin{equation}\label{eq:Def:OperatorChi}
    \OpChi(\zeta)(t):=(\zeta(t),\dot{\zeta}(t),\ldots,\zeta^{(r-1)}(t))\in\R^{rm}
\end{equation}
for a function $\zeta\in W^{r,\infty}(I,\R^m)$ on an interval $I\subset\Rp$ and $t\in I$. 
Then, we get
\begin{equation}\label{eq:ErrorsDyn}
\begin{aligned}
    \xi_1(\OpChi(\zeta)(t))&=\zeta(t),\\
    \xi_{i+1}(\OpChi(\zeta)(t))&=\dd{t}\xi_i(\OpChi(\zeta)(t))+k \xi_i(\OpChi(\zeta)(t)).
\end{aligned}
\end{equation}
\end{remark}

With the auxiliary error variable $\xi_r$ as in~\eqref{eq:ErrorVars} and $\OpChi$ 
as in~\eqref{eq:Def:OperatorChi}, we define the following novel funnel controller with constant gain terms
for a funnel function~$\psi\in\cG$ and systems of the form~\eqref{eq:Sys}:
\begin{equation}\label{eq:DefinitionFC}
\boxed{
\begin{aligned}
    e_r(t) &= \xi_r(\OpChi(y-y_{r})(t)),\\
    u(t)   &= N\circ \gamma(\|\tfrac{\alpha}{\beta}\cdot e_r(t)\|^2) e_r(t),\\
\end{aligned}
}    
\end{equation}
where $N \in \cC(\Rp,\R)$ is a surjection, 
$\gamma \in\cC([0,1) ,[1,\infty))$ is a bijection, 
and $\alpha, \beta$ are the associated constants of $\psi$.
A suitable choice for the bijection is $\gamma(s):=1\slash(1-s)$.

\begin{remark}
Although differing in certain details, the standard funnel controller designs
from~\cite{BergLe18a,BergIlch21} utilise auxiliary error variables structurally
similar to the ones in~\eqref{eq:ErrorVars}.
However, instead of a constant $k>0$, time-varying gain functions $k_1,\ldots,k_{r-1}$
of the form $k_i(t,z)=1\slash(1-\phi_i(t)^2\Norm{\xi(t,z)}^2)$ are used where
$1/\phi_i$ is an additional funnel function for the error $\xi_{i}$.
This reciprocal time-varying term introduces a non-linear dependency between the  
error signals $\xi_{i+1}$ and $\xi_i$.
Moreover, the controller~\eqref{eq:DefinitionFC} utilises the controller gain
$N\circ \gamma(\|\tfrac{\alpha}{\beta}e_r(t)\|^2)$ instead of $N\circ
\gamma(\|e_r(t)/\phi(t)\|^2)$ with a varying funnel function $\phi\in\cG$ as
in~\cite{BergLe18a,BergIlch21}.
These changes simplify the design controller design.
However, contrary to earlier works, 
the controller~\eqref{eq:DefinitionFC} cannot adjust its 
actuating signal depending on the funnel~$\psi$. 
An over some time interval widening funnel might be advantageous in some
applications, for example in the presence of periodic disturbances.
The controller~\eqref{eq:DefinitionFC} is not able to take advantage of this.
\end{remark}

We will show that  the application of the controller~\eqref{eq:DefinitionFC}
ensures that the tracking error $e$ evolves in a prescribed performance funnel
if $k>0$ is chosen sufficiently large and the initial value is sufficiently
close to the reference trajectory.

\section{Funnel control -- main result}\label{Sec:Main}
In this section, we show that the application of the funnel
controller~\eqref{eq:DefinitionFC} to system~\eqref{eq:Sys} belonging to the system
class $\cN^{r,m}$ as defined in Definition~\ref{Def:system-class} leads to a closed-loop
initial value problem which has a global solution. 
By \emph{solution} in the sense of \emph{Carath\'{e}odory} of~\eqref{eq:Sys},
we mean a function 
$x = (x_1,\ldots,x_r):[0,\omega) \to \R^{rm}$, $\omega>t_0$, where 
\[
    x|_{[0,t_0]} = (y^0,\dot y^0,\ldots,(y^0)^{(r-1)}),
\]
and $x\vert_{[t_0,\omega)}$ is absolutely continuous such that $\dot x_i(t) = x_{i+1}(t)$ for $i=1,\ldots,r-1$, and 
\[
    \dot x_r(t) = f(d(t),\oT(x(t)), u(t))
\]
for almost all~$t\in[t_0,\omega)$.
As before, the initial condition is replaced by
$x(t_0)=y^0\in\R^{rm}$ in the case $t_0=0$. 
A solution $x$ is said to be \emph{maximal}, if it has no right extension that is also a solution.
It is called \emph{global} if $\omega=\infty$.

We are now in the position to present the main result of the paper.
\begin{theorem}\label{Thm:FunCon}
Consider system~\eqref{eq:Sys} with $(d,f,\oT)\in\cN^{m,r}$ 
and initial data $y^0  \in \cC^{r-1}([0,t_0],\R^m)$.
Let $\psi\in\cG$ with associated constants $\alpha, \beta >0$ and 
choose parameter 
\begin{equation}\label{eq:def:parameters}
    k\geq \alpha+2.
\end{equation}
Furthermore, let $y_{\rf}\in W^{r,\infty}(\Rp,\R^{m})$ such that
\begin{align*}
    \Norm{\xi_1(\OpChi(y-y_{\rf})(t_0))}&<\psi(t_0),\\
    \Norm{\xi_i(\OpChi(y-y_{\rf})(t_0))}&<\beta/\alpha,\quad i=2,\ldots,r.
\end{align*}
Then, the application of the funnel controller~\eqref{eq:DefinitionFC} 
to the system~\eqref{eq:Sys} yields an initial-value problem which has a solution,
and every maximal solution $y: [0,\omega)\to\R^m$ has the following properties:
\begin{enumerate}[(i)]
    \item\label{Item:GlobalSol} The solution is global, i.e. $\omega = \infty$.
    \item\label{Item:ControlBounded} The input $u:\Rp\to\R^m$ is bounded.
    \item\label{Item:FunnelProp} The tracking error~$e=y-y_{\rf}$ evolves in its performance funnel given by~$\psi$, i.e.
    \[
        \forall\,t\geq t_0:\quad  \Norm{e(t)}<\psi(t).
    \]
\end{enumerate}
\end{theorem}
\begin{remark}
Theorem~\ref{Thm:FunCon} poses boundaries on the initial values of the tracking error to
ensure the proper operation of the funnel controller~\eqref{eq:DefinitionFC}.
Practically, this requirement can be satisfied by selecting the reference
trajectory appropriately. In contrast, the controller frameworks in
\cite{BergLe18a,BergIlch21} also stipulate that the initial values of their
auxiliary error variables reside within prescribed funnel boundaries $\phi_i$
Crucially, however, these boundaries are designer-configurable, granting
flexibility to accommodate a wider range of initial conditions.
\end{remark}
\section{Proof of the main result}
Throughout this section, let the assumptions of Theorem~\ref{Thm:FunCon} hold.
In particular, we presume that the function 
$\psi\in\cG$ with associated constants $\alpha, \beta >0$ is fixed,  
the parameter $k\geq\alpha+2$, 
and that the auxiliary error functions $\xi_i$ for $i=1,\ldots, r$ are defined as in~\eqref{eq:ErrorVars}.
We show that the application of the funnel controller~\eqref{eq:DefinitionFC} to the
system~\eqref{eq:Sys} achieves that the tracking error $e=y-y_{\rf}$ evolves within
the funnel $\cF_{\psi}$ by ensuring that $\OpChi(e)(t)$ is an element of the set 
\begin{equation}\label{eq:DefSetD}
    \cD_{t}:=\setdef{z\in\R^{rm}}{
    \begin{array}{l}
         \Norm{\xi_1(z)}<\psi(t)  \\
         \Norm{\xi_i(z)}<\tfrac{\beta}{\alpha}, i = 2,\ldots,r
    \end{array}
    }
\end{equation}
for all $t\geq t_0$.  By assumption we have $\OpChi(e)(t_0)\in\cD_{t_0}$.

Before proving Theorem~\ref{Thm:FunCon}, we show in the following lemma that it
is sufficient to ensure $\Norm{\xi_r(\OpChi(e)(t))}<\tfrac{\beta}{\alpha}$ in
order to guarantee $\OpChi(e)(t)$ is an element of $\cD_{t}$ for $t\geq t_0$.

\begin{lemma}\label{Lemma:Errors}
    Let $t_0\geq 0$, $s>t_0$, and $\zeta\in\cC^{r-1}([t_0,\infty),\R^m)$ with $\OpChi(\zeta)(t_0)\in\cD_{t_0}$.
    If $\Norm{\xi_{r}(\OpChi(\zeta))(t)}<\tfrac{\alpha}{\beta}$ for all $t\in[t_0,s)$, then
    $\OpChi(\zeta)(t)\in\cD_{t}$ for all $t\in [t_0,s)$.
\end{lemma}
\begin{pf}
    We modify the proof of~\cite[Lemma 3.2]{BergDenn24} to the current setting.
    Seeking a contradiction, we assume that there exists $t \in (t_0, s)$ such that either
    $ \| \xi_1(\OpChi(\zeta)(t)) \| \geq \psi(t)$ or $ \| \xi_i(\OpChi(\zeta)(t)) \| \geq \beta/\alpha$
    for at least one ${i \in \{2,\ldots,r-1\}}$.
    W.l.o.g.\ let~$i$ be the largest index with this property.
    We use the shorthand notation $e_i(t):=\xi_i(\OpChi(\zeta)(t))$.    
    We separately consider the two cases $i=1$ and $i>1$.
    First, we suppose $i=1$.
    Invoking continuity of~$e_1$, there exist
    $t^\star := \min \setdef{ t \in [t_0, s] }{ \| e_1(t) \| = \psi(t) }$
    and 
    $t_\star := \min\setdef{\! t \in [t_0, t^\star) }{ \fa s\in[t,t^\star]: \sqrt{2} \| e_1(s) \| \geq \psi(s)\! }$.
    Then, note that $\psi(t)\geq \Norm{e_1(t)}\geq \psi(t)/\sqrt{2}$ for all $t\in[t_\star,t^\star]$.
    By properties of~$\psi\in \cG$, we have
    \[
        -\frac{\dot{\psi}(t)}{\psi(t)}\le \frac{\alpha\psi(t) - \beta}{\psi(t)} \le \alpha,
    \]
    and
    $\psi(t) \ge \psi(t_0) \me^{-\alpha(t-t_0)} + \frac{\beta}{\alpha}>\frac{\beta}{\alpha}$ for all $t\ge t_0$.
    With~\eqref{eq:ErrorsDyn}, we have $\dot{e}_1(t)=e_2(t)-k e_1(t)$ for $t \in [t_{\star},t^\star]$.
    Thus, omitting the dependency on $t$, we calculate, for $t \in [t_{\star},t^\star]$,
    \begin{align*}
       \tfrac{1}{2}\dd{t}\Norm{\tfrac{e_1}{\psi}}^2
        &= \al\tfrac{e_1}{\psi},\tfrac{\dot{e}_1\psi-e_1\dot{\psi}}{\psi^2}\ar
        = \al\tfrac{e_1}{\psi},- \rbl k+\tfrac{\dot{\psi}}{\psi}\rbr  \tfrac{e_{1}}{\psi}+\tfrac{e_{2}}{\psi}\ar\\
        &\le  -\rbl k+\tfrac{\dot{\psi}}{\psi}\rbr\Norm{\tfrac{e_1}{\psi}}^2  + \tfrac{ \| e_1 \| \|e_2\|}{\psi^2}\\
        &\leq -\rbl k+\tfrac{\dot{\psi}}{\psi}\rbr \tfrac{1}{2} + \tfrac{\alpha}{\beta}\|e_2\| \\
        &\leq -\rbl k-\alpha\rbr \tfrac{1}{2} +1\leq 0,
    \end{align*}
    where we used $k\geq\alpha+2$ and $\Norm{e_{2}(t)}\leq \frac{\beta}{\alpha}$ for all~$t \in [t_\star,t^\star]$.
    Thus, upon integration, the contradiction
    \begin{equation*}
        \psi(t^\star) = \| e_1(t^\star) \| \le \| e_1(t_\star)\| < \psi(t_\star)
    \end{equation*}
    arises.
    Now, we consider the  case $\| e_i(t) \|\ge \beta/\alpha$ for $i>1$.
    Similar as before, define 
    $t^\star := \min \setdef{\!\! t \in [t_0, s]\!\!\! }{\!\! \| e_i(t) \| = \tfrac{\beta}{\alpha}\!\!\!}$
    and 
    $t_\star := \min\setdef{\!\! t \in [t_0, t^\star)\!\! }{\!\! \fa s\in[t,t^\star]: \sqrt{2} \| e_i(s) \| \geq \tfrac{\beta}{\alpha} \!\!}$.
    Then,  $\tfrac{\beta}{\alpha}\geq \Norm{e_i(t)}\geq \beta/(\alpha\sqrt{2})$ for all $t\in[t_\star,t^\star]$.
    Utilising~\eqref{eq:ErrorsDyn}, we calculate, for $t \in [t_{\star},t^\star]$,
    \begin{align*}
        \tfrac{1}{2 }\dd{t} \| e_i \|^2
        &=  \langle e_i , \dot{e}_i \rangle
        =  \langle e_i, -k e_i + e_{i+1} \rangle\\
        &\le - k \| e_i \|^2 + \| e_i \| \| e_{i+1} \|
        \leq \tfrac{\beta^2}{ \alpha^2} \left( - \tfrac{k}{2} + 1 \right)  \le  0,
    \end{align*}
    where we used~$k\geq 2$ and that $\| e_{i+1}(t)\| \leq \tfrac{\beta}{\alpha}$ for all~$t \in [t_\star,t^\star]$ 
    by maximality of $i$.
    Hence, the contradiction
    \[
      \tfrac{\beta^2}{\alpha^2} \leq \| e_i(t^\star)\|^2 \leq \| e_i(t_\star)\|^2< \tfrac{\beta^2}{\alpha^2}
    \]
    arises, which completes the proof.
    \xqed
\end{pf}

\noindent
\textbf{Proof of Theorem~\ref{Thm:FunCon}.}
We modify the proof of \cite[Thm.~1.9]{BergIlch21} to the current setting.

\noindent
\emph{Step 1}:
We show that, when applying the control law~\eqref{eq:DefinitionFC} to the
system~\eqref{eq:Sys} there exists a maximal solution $x=(x_1,\ldots,x_r):[0,\omega)\to\R^{rm}$, $\omega>t_0$.
Since the controller~\eqref{eq:DefinitionFC} introduces a pole on the right hand side of
the closed-loop differential equation, some care must be exercised. In order to
do that we reformulate the system in the form of an initial-value problem to
which  the well-known existence theory applies.
Define the non-empty open set
\[
    \cE:=\setdef
    {\! (t,z)\in\Rp\times\R^{rm}\!}
    {
          \Norm{\xi_{r}(z-\OpChi(y_{\rf}(t))}<\tfrac{\beta}{\alpha}
    }
\] 
and the function $F:\cE\times\R^q\to\R^{rm}$ mapping $(t,z,\eta)=(t,z_1,\ldots,z_r,\eta)$ on 
\[
   \begin{bmatrix}
        z_2\\
        \vdots\\
        z_r\\
        f(d(t),\eta,(N\circ\gamma)(\Norm{z-\OpChi(y_{\rf})}^2) \xi_r(z-\OpChi(y_{\rf})(t)))
    \end{bmatrix}. 
\]
Utilising the notation $x(t)=\OpChi(y)(t)$, the initial value problem~\eqref{eq:Sys} with 
with feedback controller~\eqref{eq:DefinitionFC} takes the form 
\begin{equation}\label{eq:InitalvalProb}
\begin{aligned}
    \dot{x}&=F(t,x(t),\oT(x)(t)),\\
     x|_{[0,t_0)}&=\OpChi(y^0)\in\cC([0,t_0],\R^{rm}).
\end{aligned}
\end{equation}
We have $(t_0,x(t_0))\in\cE$ by assumption.
Application of a variant of~\cite[Theorem~B.1]{Ilchmann01102009} yields the existence of a maximal solution
$x=(x_1,\ldots,x_r):[0,\omega)\to\R^{rm}$, $\omega>t_0$, of~\eqref{eq:InitalvalProb}
with $\graph(x|_{[t_0,\omega)})\subset\cE$.
Moreover, the closure of $\graph\rbl x|_{[t_0,\omega)}\rbr$ is not a compact subset of~$\cE$.
Note that the first $m$-dimensional component $x_1$ of $x$
is a maximal solution of the initial value problem~\eqref{eq:Sys} with control law~\eqref{eq:DefinitionFC}.
We, therefore, will denote it in the following with $y$, i.e. $y=x_1$.
Moreover, we will use the short-hand notation 
\[
    e_i(t):=\xi_i(x-\OpChi(y_{\rf})(t))
\]
for $i=1,\ldots,r$ in the following. Utilising this notation, we have
$e_1(t)=\xi_i(x-\OpChi(y_{\rf})(t))= y(t)-y_{\rf}(t)$, i.e. $e_1=e$ is the tracking error.

\noindent
\emph{Step 2}: We define some constants for later use.
As the function $F$ is defined on the set $\cE$, we have $\Norm{e_r(t)}<\beta/\alpha$ for all $t\in [t_0,\omega)$.
Since, in addition, $\OpChi(e)(t_0)\in\cD_{t_0}$, Lemma~\ref{Lemma:Errors} yields
$\Norm{e_1(t)}<\psi(t)$ and $\Norm{e_i(t)}<\beta/\alpha$ for all $t\in [t_0,\omega)$ and all $i=2,\ldots, r-1$.
Due to the definition of the error variables~$\xi_i$ in~\eqref{eq:ErrorVars} there exists an invertible matrix $S\in\R^{rm\times rm}$ such that
\begin{equation}\label{eq:reps_ei_chi}
    \begin{pmatrix} \xi_1(x- \OpChi(y_{\rf})) \\ \vdots \\  \xi_r(x- \OpChi(y_{\rf}))\end{pmatrix} = S \rbl x-\OpChi(y_{\rf})\rbr.
\end{equation}
Hence, by boundedness of $\psi$ and $y_{\rf}^{(i)}$ for all~$i=1,\ldots, r$,
there exists a compact set $K\subset\R^{rm}$ with 
\[
    \fa t\in[t_0,\omega):\quad x(t)\in K.
\]
Invoking the BIBO property of operator~$\oT$, there exists a compact set $K_q\subset\R^q$ with  
$\oT(x)(t)\in  K_q$ for all $t\in [t_0,\omega)$.
Further, there exists a compact set $K_p\subset\R^p$ with $d(t)\in K_p$ for almost all $t\in[t_0,\omega)$
since  $d\in L^\infty([t_0,\infty),\R^p)$. 
There exists $\nu\in(0,1)$ such that the function
$\mathfrak{h}:\R\to\R$ given by 
\[
\mathfrak{h}(s):=\min\setdef{\langle v, f(\delta,z, -s v)\rangle \!\!}
{\!\!\! 
\begin{array}{ll}
     \delta\in K_p,  z \in  K_q,   \\
v\in\R^m,\nu \leq \|v\| \leq 1
\end{array}
\!\!\!}
\]
is unbounded from above because of the high-gain property of the function $f$, see Definition~\ref{Def:system-class}.
Define 
$\mu_1^0 := \SNorm{\psi}$ and $\mu_i^0 := \beta/\alpha$ for $i=2,\ldots, r$.
Moreover, 
\begin{equation}\label{eq:DefMus}
    \mu_{i}^{j+1}:= \mu_{i+1}^{j}+k\mu_{i}^{j}
\end{equation}
for  $i=1,\ldots, r$ and $j=0,\ldots,r-i-1$.
Due to the unboundedness of the function~$\mathfrak{h}$ and the surjectivity of~$N\circ\gamma$ it is possible to choose $\eps\in(0,1)$ 
such that $\eps>\max\cbl{\nu,\| \tfrac{\alpha}{\beta} e_r(t_0)\|}\cbr$, and 
\[
    \tfrac{1}{2}\mathfrak{h}((N\circ\gamma)(\eps^2))\geq\lambda:=\SNorm{y_{\rf}^{(r)}}+\sum_{j=1}^{r-1}k\mu_{j}^{r-j}.
\]

\noindent
\emph{Step 3}:  We show $\|\tfrac{\alpha}{\beta} e_r(t)\|\leq\eps$ for all $t\in [t_0,\omega)$. 
Seeking a contradiction, assume there exists $t^{\star}\in [t_0,\omega)$ with $\Norm{e_r(t^\star)}>\eps$.
As $\|\tfrac{\alpha}{\beta}e_r(t_0)\|<\eps$, there exists 
\[
    t_{\star}:=\sup\setdef{t\in[t_0,t^\star)}{\|\tfrac{\alpha}{\beta}e_r(t)\|=\varepsilon}<t^{\star}
\]
because $e_r$ is a continuous function on $[t_0,t^\star]$.
By definition of $t_\star$, we have 
$\|\tfrac{\alpha}{\beta}e_r(t)\|\geq\eps>\nu$ for all $t\in [t_\star,t^\star]$ and 
\[
    \mathfrak{h}((N\circ\gamma)(\|\tfrac{\alpha}{\beta}e_r(t_\star)\|^2))\geq2\lambda.
\]
Thus, there exists $\hat{t}\in[t_\star,t^\star]$ with 
\[
    \fa t\in[t_\star,\hat{t}]:\quad\mathfrak{h}((N\circ\gamma)(\|\tfrac{\alpha}{\beta}\,e_r(t)\|^2))\geq\lambda.
\]
Note that 
$\Norm{e_1(t)}<\psi(t)\leq \mu_1^0$ and 
$\Norm{e_i(t)}<\beta/\alpha= \mu_i^0$
on the interval $[t_0,t^{\star})$ for all $i=2,\ldots, r$.
Using the definition of the error variable $\xi_i$ in~\eqref{eq:ErrorVars} and the definition of~$\mu_i^j$ it follows that
\begin{equation*}
    \Norm{e^{(j+1)}_i(t)}
    \!=\!
    \Norm{e^{(j)}_{i+1}(t)-ke^{(j)}_{i}(t)}
    \!\leq
    \mu_{i+1}^j\!+k\mu^{j}_{i}
    \!=\!\mu^{j+1}_{i}
\end{equation*}
inductively for all $i=1,\ldots, r$ and $j=0,\ldots,r-i-1$.
Moreover, utilising the definition of the error variable $e_r$ in~\eqref{eq:ErrorVars}, one can show by induction that
\begin{equation*}
        e_r(t)=e^{(r-1)}(t)+\sum_{j=1}^{r-1}ke_j^{(r-j-1)}(t).
\end{equation*}
Omitting the dependency on $t$, we calculate
\begin{align*}
&\dd{t} \tfrac{1}{2} \Norm{ e_r}^2
= \al e_r, \dot{e}_r\ar
 = \al e_r, e^{(r)}+\sum_{j=1}^{r-1}ke_j^{(r-j)}\ar\\
& =  \al e_r, f(d,\oT(x),u)-y_{\rf}^{(r)}+\sum_{j=1}^{r-1}ke_j^{(r-j)}\ar\\
& <   \| y^{(r)}_{\rf}\|_{\infty} +\sum_{j=1}^{r-1}k\mu_{j}^{r-j}+ \al  e_r , f(d, \oT(x), u) \ar  \\
& =  \lambda 
+ 
\al  e_r, f\rbl d, \oT(x),  (N\circ \gamma)\rbl\|\tfrac{\alpha}{\beta} e_r\|^2\rbr e_r \rbr \ar\\
& \le \lambda
-  \min \setdef{\!\!\!\! \al v, f\rbl\delta,z, - (N \! \circ \! \gamma)\!\rbl\|\tfrac{\alpha}{\beta} e_r\|^2\rbr v\rbr \!\ar \!\!\!}
{  \!\! \begin{array}{l}
     \delta \!\in\! K_{{p}},  \\
     z \!\in\! K_q, \\
     v \!\in\! V
\end{array}\!\!\!\!} \\
& \le \lambda  -  \mathfrak{h}\rbl(N \circ \gamma)\rbl\|\tfrac{\alpha}{\beta} e_r\|^2\rbr\rbr
 \le 0,
\end{align*}
for almost all $t \in [t_\star,\hat{t}]$.
Integration yields 
$
    \eps < \|\tfrac{\alpha}{\beta} e_r(\hat{t})\| \le \|\tfrac{\alpha}{\beta} e_r(t_\star)\| = \eps
$,
a contradiction.
Therefore, $\|\tfrac{\alpha}{\beta} e_r(t)\|\leq\eps$ for all $t\in [t_0,t^{\star}]$.
This is a contradiction to the definition of $t^{\star}$.
Thus, $\|\tfrac{\alpha}{\beta}e_r(t)\|\leq\eps$ for all $t\in [t_0,\omega)$.

\noindent
\emph{Step 4}:
As a consequence of and Step~3 and Lemma~\ref{Lemma:Errors}, we have 
$\|\tfrac{\alpha}{\beta}e_r(t)\|\leq\eps$ and $\OpChi(e)(t)\in\cD_t$ for all $t\in [t_0,\omega)$.
By the definition of~$\xi_i$ in~\eqref{eq:ErrorVars} for $i=1,\ldots, r$ and the boundedness of the
function~$\OpChi(y_{\rf})$, the solution $x$ is a bounded function, too.
Since the closure of $\graph\rbl x|_{[t_0,\omega)}\rbr$ is not a compact subset of $\cE$, 
this implies $\omega=\infty$ and thereby shows assertion~\eqref{Item:GlobalSol}.
Further, $\|\tfrac{\alpha}{\beta}e_r(t)\|\leq \eps<1$ for all $t\in[t_0,\infty)$ implies the boundedness of $u$ in~\eqref{eq:DefinitionFC},
i.e. assertion~\eqref{Item:ControlBounded}.
Moreover, by $\OpChi(e)(t)\in\cD_t$ for all $t\in [t_0,\infty)$, we have
$
    \Norm{y(t)-y_{\rf}(t)}=\Norm{\xi_1(\OpChi(y-y_{\rf}))(t)}<\psi(t)
$ for all $t\in[t_0,\infty)$.
This shows assertion~\eqref{Item:FunnelProp} and completes the proof.
\hfill$\Box$

\section{Simulations}
\begin{figure}[htp]
    \begin{center}
    \includegraphics[trim=2cm 4cm 5cm 15cm,clip=true,width=5.5cm]{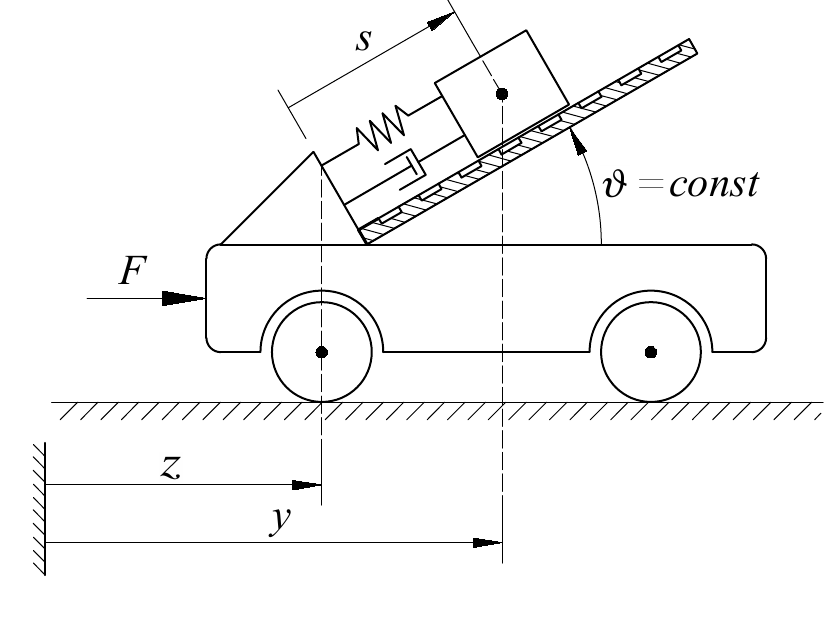}
    \end{center}
    \vspace{-4mm}
    \caption{Mass-on-car system.The figure is based on~\cite{SeifBlaj13}}
    \label{Fig:Mass_on_car}
\end{figure}

We compare the the funnel controller~\eqref{eq:DefinitionFC} to the controller designs presented in~\cite{BergLe18a} 
and \cite{BergIlch21}. 
To this end, we consider the example of a mass-spring system mounted on a car from~\cite{SeifBlaj13}, 
which has already been used as an example in the earlier works.
As depicted in Figure~\ref{Fig:Mass_on_car}, a mass~$m_2$ 
is coupled to a car with mass~$m_1$ via a spring-damper component and 
moves frictionless on a ramp inclined by the angle~$\vartheta$.
The coupling component is parametrised by the spring constant $c>0$ and the damping coefficient $\delta>0$. 
It is assumed that force $u= F$ can be applied as a control to the car.
The dynamics of the system can be described by the equations
\begin{small}
\begin{align}\label{eq:ExampleMassOnCarSystem}
    \begin{bmatrix}
        m_1 + m_2& m_2\cos(\vartheta)\\
        m_2 \cos(\vartheta) & m_2
    \end{bmatrix}\!\!
    \begin{pmatrix}
        \ddot{z}(t)\\
        \ddot{s}(t)
    \end{pmatrix}
    \!+\!
    \begin{pmatrix}
        0\\
        c s(t) +d\dot{s}(t)
    \end{pmatrix}
    \!=\!
    \begin{pmatrix}
        u(t)\\
        0
    \end{pmatrix}\!.
\end{align}
\end{small}
Here $z(t)$ is the horizontal position of the car and $s(t)$ the relative position of the mass on the ramp at time~$t$. 
The horizontal position of the mass on the ramp, i.e.,
\[
    y(t) = z(t) +s(t)\cos(\vartheta),
\]
is assumed to be the output of the system.
As outlined in~\cite[Section 3.1]{BergIlch21}, the system belongs to the
class~$N^{1,r}$ with $r=2$ for $\vartheta\in(0,\pi/2)$ and $r=3$ for
$\vartheta=0$. As in~\cite[Section 3.1 Case 2]{BergIlch21}, we choose the
parameters $m_1=4$, $m_2=1$, $c=2$, $\delta=1$, $\vartheta=0$.
The control objective is to track the reference signal $y_{\rf}(t)=\cos(t)$
within predefined boundaries given by the  function~$\psi\in\cG$. This means
guaranteeing that the tracking error~$e(t)=y(t)-y_{\rf}(t)$ fulfils
$\Norm{e(t)}<\psi(t)$ for all $t\geq0$. We choose the initial values 
$z(0)=-0.3$, $\dot{z}=-0.21$, and $s(0)=\dot{s}(0)=1$.
As funnel function, we choose $\psi:\Rp\to\R$, $\psi(t)=3\me^{-t}+0.1$ with  
associated constants $\alpha=1$ and $\beta=0.1$.
Therefore, the gain $k=3$ is chosen for the controller~\eqref{eq:DefinitionFC}.
\begin{figure}[ht] \centering
    \begin{subfigure}[b]{0.48\textwidth}
     \centering
        \includegraphics[width=8.7cm]{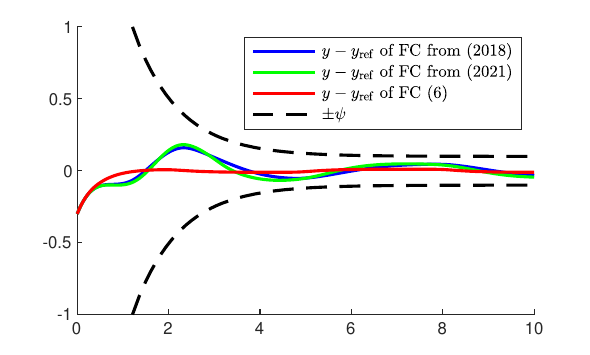}
        \caption{Tracking error~$e$ and funnel boundary~$\psi$}
        \label{Fig:SimulationOutputError}
    \end{subfigure}
    \begin{subfigure}[b]{0.48\textwidth}
        \centering
        \includegraphics[width=8.7cm]{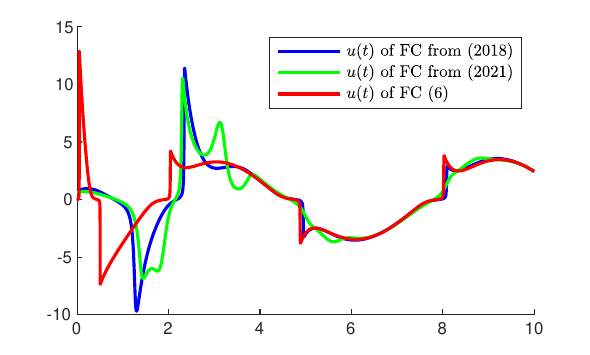}
        \caption{Control input}
        \label{Fig:SimulationControlInput}
    \end{subfigure}
    \caption{\label{Fig:SimulationFC}Simulation of system~\eqref{eq:ExampleMassOnCarSystem} under control law~\eqref{eq:DefinitionFC} and funnel control from
    \cite{BergLe18a}  and \cite{BergIlch21}.}
\end{figure}
All simulations are performed with \textsc{Matlab} on 
the interval~$[0,10]$ and depicted in Figure~\ref{Fig:SimulationFC}.
The tracking errors {resulting from the application of funnel controller~\eqref{eq:DefinitionFC} 
and both the controller designs from \cite{BergLe18a} and \cite{BergIlch21} are shown in Figure~\ref{Fig:SimulationOutputError}.
The corresponding control signals are displayed in Figure~\ref{Fig:SimulationControlInput}.
It is evident that all three control schemes achieve the control objective,
the evolution of the tracking error with in the performance boundaries given by~$\psi$.
The controller~\eqref{eq:DefinitionFC} keeps the system output
over the whole time interval close to the reference signal. The other two controllers, however, let
the error signal evolve more freely within the funnel boundaries.
This results in higher control signals when the tracking error gets close the boundary.
The error signals of the controller~\eqref{eq:DefinitionFC} are amplified by the
constant gain $k$, which for small deviations from the reference signal is large
in comparison to the the other two controllers' dynamically computed gains. 
The controller therefore counteracts deviations from $y_{\rf}$ more aggressively 
resulting in a better tracking tracking performance.
On the other hand, as Theorem~\ref{Thm:FunCon} requires $\OpChi(y-y_{\rf})(t_0)\in\cD_{t_0}$,
the initial values have to be already very close to the reference signal to guarantee 
the functioning of controller~\eqref{eq:DefinitionFC}. 
The other two controllers allow more leeway for the initial values.

\section{Conclusion}
We proposed output derivative feedback controller that achieves reference
tracking with prescribed performance for higher order systems. By utilising
auxiliary error signals only depend linear on the tracking error and its
derivatives it simplifies previously proposed controllers
designs~\cite{BergLe18a,BergIlch21}. This has the potential of reducing
numerical difficulties as the usage of reciprocal intermediate penality terms is
avoided. However, the controller loses its ability to dynamically react to
changing funnel boundaries. Future research will focus on how to recover this
property while keeping the simplified controller structure.

\section{Acknowledgements} 
I thank Thomas Berger (University of Paderborn) and Lukas Lanza (TU Ilmenau) for many fruitful discussions.
\small
\bibliography{references}

\end{document}